\documentclass[12pt]{article}
\usepackage[latin1]{inputenc}
\usepackage{amsmath,amsthm,amssymb}

\newtheorem{thm}{Theorem}
\numberwithin{thm}{section}

\def\eq#1{(\ref{#1})}

\newcommand{\neweq}[1]{\begin{equation}\label{#1}}
\def\weak{\rightharpoonup}

\def\phi{\varphi}
\def\RR{\mathbb R}

\def\di{\displaystyle}
\def\ri{\rightarrow}
\def\intom{\int_\Omega}
\def\huo{H^1_0(\Omega )}
\def\incep{\left\{\begin{array}{cl} }
 \def\termin{\end{array}\right. }
\def\2af{2^*_\alpha}

\def\proof{{\it Proof.}\ }

\title{\sc Existence results for variational-hemivariational
problems with lack of convexity}
\author{\sc Vicen\c tiu R\u adulescu$\,^{a,b}$ and Du\u san Repov\v s$\,^{c}$\\
\small $^a\,$Institute of Mathematics ``Simion Stoilow" of the Romanian
Academy, Bucharest\\
\small $^b\,$Department of Mathematics,
University of Craiova,  200585 Craiova, Romania\\
 \small $^c\,$Institute of Mathematics,
Physics and Mechanics, University of Ljubljana,\\ \small Jadranska  19,  P. O. Box 2964, 1001 Ljubljana, Slovenia\\
 \small E-mail: {\tt vicentiu.radulescu@imar.ro}\qquad {\tt dusan.repovs@guest.arnes.si}\\}

\date{}

\begin{document}

\maketitle
\begin{abstract}
We establish existence results of Hartmann-Stampacchia type for a class of variational-hemivariational
inequalities on closed and convex sets (either bounded or unbounded) in a Hilbert space.
\end{abstract}

\section{Introduction}

Let $\Omega$ be a bounded open set in $\RR^N$. Assume that $K$ is a nonempty, closed, and convex (bounded or unbounded) set in $H^1_0(\Omega)$. The first major result in the theory of variational inequalities is the following direct consequence of the Stampacchia theorem: for any $f\in H^{-1}(\Omega)$, there is a unique $u\in K$ such that for  all $v\in K$,
\begin{equation}\label{stampa}
\int_\Omega \nabla u\cdot\nabla (v-u)\,dx\geq \langle f,v-u\rangle\,.\end{equation}
The above result is often referred as the Hartman-Stampacchia theorem (see \cite[Lemma~3.1]{ref6} or \cite[Theorem I.3.1]{ref9}). A simple proof of the Hartmann-Stampacchia theorem is due to Brezis and may be found in \cite{ref9}.

Several nonlinear and nonconvex extensions of \eq{stampa} have been given in a nonsmooth framework by Fundos, Panagiotopoulos and R\u adulescu \cite{fpr} and by Motreanu and R\u adulescu \cite{mrnfa}. We refer to \cite{aga1}, \cite{aga2}, \cite{krv}, \cite{mopa} for related results and applications.

In \cite{fpr} there are obtained Hartman-Stampacchia type properties for nonconvex inequality problems of the type: find $u\in K$ such that for all $v\in K$
$$\int_\Omega \nabla u\cdot\nabla (v-u)\,dx+\int_\Omega j^0(x,u(x);v(x)-u(x))\,dx\geq 0,$$
where $j^0$ stands for the Clarke generalized directional derivative. The  case of variational-hemivariational inequalities was studied in \cite{ref6} for the model problem: find $u\in K$ such that for all $v\in K$
$$\int_\Omega \nabla u\cdot\nabla (v-u)\,dx+\Phi(v)-\Phi(u)+\int_\Omega j^0(x,u(x);v(x)-u(x))\,dx\geq 0,$$
where $\Phi$ is convex and lower semicontinuous.

In the present paper we are concerned with a more general class of inequality problems with lack of convexity. The main idea in the study we develop in this work is related with the previous nonlinear inequality problems but is also in strong relationship with the semilinear boundary value problem
\begin{equation}\label{ee1}
\left\{\begin{array}{lll}
&\di -\Delta u=f(x,u)&\qquad\mbox{in $\Omega$}\\
&\di u=0&\qquad\mbox{on $\partial\Omega$},\end{array}\right.
\end{equation}
where $f:\Omega\times\RR\ri\RR$ is a Carath\'eodory function satisfying
\begin{equation}\label{ee2}
|f(x,t)|\leq\lambda_1\,|t|\qquad\mbox{for all $(x,t)\in\Omega\times\RR$.}
\end{equation}
Here, $\lambda_1$ denotes the first eigenvalue of the Laplace operator $(-\Delta)$ in $H^1_0(\Omega)$. If $\varphi_1$ is a positive eigenfunction of $(-\Delta)$ corresponding to $\lambda_1$ then, by our basic assumption \eq{ee2}, $\varphi_1$ (resp, $-\varphi_1$) is a super-solution (resp., a sub-solution) of problem \eq{ee1}. Thus, problem \eq{ee1} has at least one solution. However, we point out that assumption \eq{ee2} is very sensitive, in the sense that problem \eq{ee1} has no longer solutions provided that $f$ has a growth described by $|f(x,t)|\leq\lambda_1\,|t|+C$, for some $C>0$. For instance, the linear Dirichlet problem
$$\left\{\begin{array}{lll}
&\di -\Delta u=\lambda_1u+1&\qquad\mbox{in $\Omega$}\\
&\di u=0&\qquad\mbox{on $\partial\Omega$}\end{array}\right.$$
does not have any solution, as can be easily seen after multiplication with $\varphi_1$.

We intend to show in the present paper that the growth assumption \eq{ee2} can be used to obtain existence results for a general class of variational-hemivariational inequalities.

\section{The main result}
We first recall that if $\varphi :H^1_0(\Omega)\ri\RR$
is a locally Lipschitz function then $\varphi^0(u;v)$ denotes the Clarke generalized derivative of $\varphi$ at $u\in H^1_0(\Omega)$ with respect to
the direction $v\in H^1_0(\Omega)$, that is, $$ \varphi^{0}(u;v) =
\limsup\limits_{w \to u \atop \lambda\downarrow 0}\frac{\varphi (w
+ \lambda v)-\varphi (w)}{\lambda}\, . $$ Accordingly, Clarke's
generalized gradient $
\partial \varphi (u)$ of $\varphi$ at $u$ is defined by $$
\partial \varphi (u) = \lbrace \xi \in H^{-1}(\Omega) \, ; \ \langle \xi
, v \rangle \leq \varphi^{0}(u;v), \  \mbox{for all}\ v \in H^1_0(\Omega) \rbrace. $$

The function $(u,v)\longmapsto \varphi^0(u,v)$ is upper semicontinuous and
$$\varphi^0(u;v)=\max\{\langle\zeta ,v\rangle;\ \zeta\in\partial\varphi (u)\}\qquad\mbox{for all $v\in H^1_0(\Omega)$}.$$
Then $\partial\varphi (u)$ is a nonempty, convex, and weak $*$ compact subset of $H^{-1}(\Omega)$.

We refer to the monograph Clarke \cite{clarke} for further properties of the generalized gradient of locally Lipschitz functionals.

\medskip
In this paper we are concerned with the following inequality problem:
\begin{equation}\label{p}
\left\{\begin{array}{ll}
&\di\mbox{find $u\in K$ such that for all $v\in K$,}\\
&\di \int_\Omega \nabla u\cdot\nabla (v-u)\,dx+\int_\Omega f(x,u)(v-u)\,dx+
\int_\Omega j^0(x,u(x);v(x)-u(x))\,dx\geq 0.
\end{array}\right.
\end{equation}

Throughout we assume that $f:\Omega\times\RR\ri\RR$ is a Carath\'eodory function such that
\begin{equation}\label{f1}
\sup_{x\in\Omega}\limsup_{t\ri\pm\infty}\left|\frac{f(x,t)}{t}\right|<\lambda_1.\end{equation}

Observe that assumption \eq{f1} implies the existence of some $\mu\in (0,\lambda_1)$ and $C>0$ such that for all $(x,t)\in\Omega\times\RR$,
\begin{equation}\label{f2}
|f(x,t)|\leq \mu\,|t|+C.\end{equation}

We assume that $j:\Omega\times\RR\ri\RR$ is a Carath\'eodory function such that
\begin{equation}\label{j1}
|j(x,y_1)-j(x,y_2)|\leq k(x)\,|y_1-y_2|\qquad\mbox{for all $x\in\Omega$ and $y_1,\, y_2\in\RR$,}
\end{equation}
for some function $k\in L^2(\Omega)$,
and there exist $h_1\in L^2(\Omega)$ and $h_2\in L^\infty (\Omega)$ such that
\begin{equation}\label{j2}
|z|\leq h_1(x)+h_2(x)|y|\qquad\mbox{for all $(x,y)\in\Omega\times\RR$ and all $z\in\partial j(x,y)$.}
\end{equation}

Our main result in this paper is the following.

\begin{thm}\label{th1}
Assume that $K$ is a nonempty, closed, and convex  set in $H^1_0(\Omega)$ and that hypotheses \eq{f1}, \eq{j1} and \eq{j2} are fulfilled. Then problem \eq{p} has at least one solution.
\end{thm}

We conclude this section by observing that condition \eq{f1} is very related to the growth assumption \eq{ee2}. However, due to the presence in \eq{p} of the nonconvex term $\int_\Omega j^0(x,u(x);v(x)-u(x))\,dx$, we are not able to work under the same hypothesis, that is,
\begin{equation}\label{assno}
\sup_{x\in\Omega}\limsup_{t\ri\pm\infty}\left|\frac{f(x,t)}{t}\right|
\leq\lambda_1.\end{equation}
However, the techniques we use in what follows enable us to obtain the same result as stated in Theorem \ref{th1} provided that \eq{assno} holds, but
$$|f(x,t)|\leq\mu\, |t|+C\qquad\mbox{for all $(x,t)\in\omega\times\RR$},$$
for some $\mu\in (0,\lambda_1)$, where $\omega\subset\Omega$ and $|\omega|>0$.

\section{An auxiliary result}
Throughout this section we assume that $\Omega$ is bounded and we prove that the existence result stated in Theorem \ref{th1} is valid in this particular case.

Let $J:L^2(\Omega)\ri\RR$ be the mapping defined by $J(u)=\int_\Omega j(x,u(x))dx$. Our assumption \eq{j2} implies that $J$ is locally Lipschitz on $L^2(\Omega)$ and for all $u,\, v\in L^2(\Omega)$,
\begin{equation}\label{vr1}
\int_\Omega j^0(x,u(x);v(x))dx\geq J^0(u;v).\end{equation}
Since $H^1_0(\Omega)$, we obtain that relation \eq{vr1} holds for any $u,\, v\in H^1_0(\Omega)$.

We recall (see \cite{fpr}) that, in view of our assumptions \eq{j1}, \eq{j2}, and \eq{f1}, the mapping
$$H^1_0(\Omega)\times H^1_0(\Omega)\ni (u,v)\longmapsto \int_\Omega j^0(x,u(x);v(x))dx\quad\mbox{is weakly upper semicontinuous} $$
and for all $v\in H^1_0(\Omega)$, the mapping
$$H^1_0(\Omega)\ni u\longmapsto\int_\Omega f(x,u)(v-u)dx\quad\mbox{is weakly continuous.} $$

The main result of this section is the following.

\begin{thm}\label{th2}
Assume that $K$ is a nonempty, closed, convex, and bounded  set in $H^1_0(\Omega)$ and that hypotheses \eq{f1}, \eq{j1} and \eq{j2} are fulfilled. Then problem \eq{p} has at least one solution.
\end{thm}

The proof of this existence property relies on the celebrated Knaster-Kuratowski-Mazur\-kiewicz principle. We first recall that if $E$ is a vector space then a subset $A$ of $E$ is said to be finitely closed if its intersection with any finite-dimensional linear manifold $L\subset E$ is closed in the Euclidean topology of $L$. Let $X$ be an arbitrary subspace of $E$. A multivalued mapping $G:X\ri{\mathcal P}(E)$ is called a KKM-mapping if
$$\mbox{conv}\, \left\{x_1,\ldots ,x_n\right\}\subset\bigcup_{i=1}^nG(x_i)$$
for any finite set $\left\{x_1,\ldots ,x_n\right\}\subset X$.

For the convenience of the reader we recall the KKM-principle of Knaster, Kuratowski, and Mazur\-kiewicz (see \cite{dugu} and \cite{kkm}).

\begin{thm}\label{th3}
Let $E$ be a vector space, $X$ be an arbitrary subspace of $E$, and $G:X\ri{\mathcal P}(E)$ be a KKM-mapping such that $G(w)$ is finitely closed for any $w\in X$. Then the family $\{G(w)\}_{w\in X}$ has the finite intersection property.
\end{thm}

\proof We claim that it is enough to show that the inequality problem
\begin{equation}\label{pp}
\left\{\begin{array}{ll}
&\di\mbox{find $u\in K$ such that for all $v\in K$,}\\
&\di \int_\Omega \nabla u\cdot\nabla (v-u)\,dx+\int_\Omega f(x,u)(v-u)\,dx+
J^0(u;v-u)\,dx\geq 0
\end{array}\right.
\end{equation}
has a solution. This fact combined with relation \eq{vr1} implies that problem \eq{p} has at least one solution.

Returning to problem \eq{pp}, let $G:K\ri {\mathcal P}(H^1_0(\Omega))$ be the multivalued mapping defined as follows: for any $w\in\huo$, let $G(w)$ be the set of all $v\in K$ such that
$$\intom\nabla v\cdot\nabla (w-v)dx+\intom f(x,v)(w-v)dx+J^0(v;w-v)\geq 0.$$

{\sc Step 1.} {\it The set $G(w)$ is weakly closed}.

Indeed, let us assume that $v_n\in G(w)$ and $v_n\weak v$ in $\huo$. Then
$$\intom\nabla v\cdot\nabla (v-w)dx\leq\liminf_{n\ri\infty}\intom\nabla v_n\cdot\nabla (v_n-w)dx$$
and
$$\lim_{n\ri\infty}\intom f(x,v_n)(w-v_n)dx=\intom f(x,v)(w-v)dx\,.$$
Using now the upper semi-continuity of the mapping $J^0(\cdot\, ;\,\cdot )$ we obtain
$$\limsup_{n\ri\infty}J^0(v_n;w-v_n)\leq J^0(v;w-v).$$
Therefore
$$J^0(v;w-v)\geq-\liminf_{n\ri\infty}\left(-J^0(v_n;w-v_n) \right)\,.$$

Using these relations we conclude that if $v_n\in G(w)$ and $v_n\weak v$ then
$$\int_\Omega \nabla v\cdot\nabla (w-v)\,dx+\int_\Omega f(x,v)(w-v)\,dx+
J^0(v;w-v)\,dx\geq 0,$$
which shows that $v\in G(w)$. Now, using the basic assumption that $K$ is bounded, we deduce that $G(w)$ is weakly closed.

\smallskip
{\sc Step 2.} {\it  $G$ is KKM-mapping}.

Arguing by contradiction, we find $w_1,\ldots, w_n\in K$ and $z\in\mbox{conv}\, \{w_1,\ldots ,w_n\}$ such that $z\notin\cup_{j=1}^nG(w_j)$. This means that for all $j=1,\ldots ,n$,
$$\intom\nabla z\cdot\nabla (z-w_j)dx+\intom f(x,z)(z-w_j)dx+J^0(z;w_j-z)<0.$$
This means that $w_j\in C$, where
$$C:=\left\{w\in K;\ \intom\nabla z\cdot\nabla (z-w)dx+\intom f(x,z)(z-w)dx+J^0(z;w-z)<0\right\}\,.$$
Since the mapping $J^0(u;\cdot )$ is subadditive and positive homogeneous (see \cite{clarke}), the set $C$ is convex, hence $z\in C$, a contradiction.

\smallskip
{\sc Step 3.} {\it The family $\{G(w)\}_{w\in K}$ has the finite intersection property}.

This follows by combining Step 2 with Theorem \ref{th3} of Knaster, Kuratowski, and Mazur\-kiewicz. Thus, there exists $u\in\cap_{w\in K}G(w)$ or, equivalently,
$$\intom\nabla u\cdot\nabla (v-u)dx+\intom f(x,u)(v-u)dx+J^0(u;v-u)\geq 0,$$
for all $v\in K$. This concludes the proof of Theorem \ref{th2}.\qed

\section{Proof of Theorem \ref{th1}}

We apply some ideas developed in \cite{fpr} and \cite{mrnfa} which rely essentially on Theorem \ref{th2} combined with the possibility to approximate the set $K$ with bounded sets having the same structure.

Without loss of generality we assume that $0\in K$. For any positive integer $n$, set
$$K_n:=\{ w\in K;\ \|w\|\leq n\}\,.$$
Thus, $0\in K_n$ for all $n\geq n_0$, where $n_0$ is a positive integer.

Applying Theorem \ref{th2} we find $u_n\in K_n$  ($n\geq n_0$) such that for all $v\in K_n$,
\begin{equation}\label{final}\intom\nabla u_n\cdot\nabla (v-u_n)dx+\intom f(x,u_n)(v-u_n)dx+\intom j^0(x,u_n(x);v(x)-u_n(x))\,dx\geq 0.\end{equation}

We claim that the sequence $(u_n)$ is bounded in $H^1_0(\Omega)$. Arguing by contradiction and passing eventually to a subsequence, we can assume that $\|u_n\|_{H^1_0(\Omega)}\ri\infty$ as $n\ri\infty$. Taking now $v=0$ as test function in relation \eq{final} we obtain (using also our assumption \eq{f1})
\begin{equation}\label{eq7}
 \|u_n\|_{H^1_0(\Omega)}^2\di =\intom |\nabla u_n|^2dx+\intom f(x,u_n)u_ndx\leq \left| \intom j^0(x,u_n(x);-u_n(x))\,dx\right|\,.\end{equation}

Using now condition \eq{j1} we find
\begin{equation}\label{eq8}\begin{array}{ll}
\di\left| \intom j^0(x,u_n(x);-u_n(x))\,dx\right|&\di\leq\intom k(x)\, |u_n(x)|\,dx\\
&\di\leq\|k\|_{L^2(\Omega)}\,\| u_n\|_{L^2(\Omega)}\leq C\, \|k\|_{L^2(\Omega)}\,\| u_n\|_{H^1_0(\Omega)}\,,\end{array}
\end{equation}
where $C>0$ is a constant determined by the continuous embedding $H^1_0(\Omega)\subset L^2(\Omega)$.

On the other hand, our assumption \eq{f1} implies
\begin{equation}\label{eq9}
\left|\intom f(x,u_n)u_ndx\right|\leq\mu\intom u_n^2dx+C\,|\Omega |\leq\frac{\mu}{\lambda_1}\, \|u_n\|_{H^1_0(\Omega)}^2+C\, |\Omega|\,.\end{equation}

Combining relations \eq{eq7}--\eq{eq9} we obtain
$$\left(1-\frac{\mu}{\lambda_1}\right)\,\|u_n\|_{H^1_0(\Omega)}^2-C\, |\Omega|\leq C\, \|k\|_{L^2(\Omega)}\,.$$
Since $\mu\in (0,\lambda_1)$, this relation shows that the sequence $(u_n)$ is bounded in $H^1_0(\Omega)$. Thus, up to a subsequence, $u_n\weak u\in K$ in $H^1_0(\Omega)$. To conclude the proof, it remains to show that $u$ is solution of problem \eq{p}. As we have already observed in the proof of Theorem \ref{th3}, it is enough to show that $u$ verifies \eq{pp}. Fix $v\in K$. Thus, there is a positive integer $N$ such that for all $n\geq N$, $v\in K_n$. Using now Theorem \ref{th2} we find that for all $n\geq N$,
\begin{equation}\label{101}\intom\nabla u_n\cdot\nabla (v-u_n)dx+\intom f(x,u_n)(v-u_n)dx+ J^0(u_n;v-u_n)\geq 0\,.\end{equation}
Next, since $u_n\weak u$, we obtain
\begin{equation}\label{102}
\intom f(x,u)(v-u)dx=\lim_{n\ri\infty}\intom f(x,u_n)(v-u_n)dx\,,\end{equation}
\begin{equation}\label{103}
J^0(u;v-u)\geq\limsup_{n\ri\infty}J^0 (u_n;v-u_n)\end{equation}
and
$$
\intom\nabla u\cdot\nabla (u-v)dx\leq\liminf_{n\ri\infty}\intom\nabla u_n\cdot\nabla (u_n-v)dx\,,$$
hence
\begin{equation}\label{104}
\intom\nabla u\cdot\nabla (v-u)dx\geq\limsup_{n\ri\infty}\intom\nabla u_n\cdot\nabla (v-u_n)dx\,.\end{equation}

Using now relations \eq{102}--\eq{104} and passing at ``$\limsup$" in \eq{101} we conclude that $u$ solves problem \eq{pp}, so $u$ is a solution of \eq{p}. This completes the proof of Theorem \ref{th1}.\qed

\medskip
{\bf Acknowledgments}. V. R\u adulescu has been supported by Grant
CNCSIS PNII--79/2007 {\it ``Procese Neliniare Degenerate \c si
Singulare"}. D.~Repov\v s was supported
by the Slovenian Research Agency grants P1-0292-0101 and J1-9643-0101.

\end{document}